\input amsppt.sty
\magnification=1200
\hcorrection{.25in}
\vcorrection{-.35in}

\topmatter
\title Free probability for probabilists
\endtitle
\author Philippe Biane \endauthor
\address CNRS,
DMI, \' Ecole Normale
Sup\'erieure\endgraf
45, rue d'Ulm, 75005 Paris FRANCE
\endaddress
\abstract  This  is an introduction to some of the most
probabilistic aspects of free probability theory.
 
\endabstract
\endtopmatter
\document
\head{ Introduction }\endhead 
Free probability is a non-commutative probability theory, in which the
concept of
independence of classical probability is replaced by that of freeness. This
new concept
incorporates both the probabilistic idea of non-correlation involved in the
independence of
random variables, and the algebraic notion of abscence of relations, 
 like e.g. between  the generators of a (non-abelian)
  free group.
It was introduced abstractly by D. Voiculescu around 1983, in order to
study some
problems in the theory of von Neumann algebras, but some years later, around
1990, he
realized that a concrete probabilistic model of free random variables is
afforded by large
independent random matrices. While this discovery has lead to an impressive
series of
advances in von Neumann algebras, some long standing open problems being
solved with
these new ideas, we shall not discuss these applications here, rather  the
purpose of these
lectures is to try to explain to a probabilist audience why the theory of
free random
variables is an interesting  and beautiful subject in itself, with deep
connections with
classical probability, and other related areas, such as harmonic analysis,
and combinatorics. We hope that this brief introduction will help the reader
find its way into
 the litterature on free probability.
We shall  take as a departure point a natural problem
about
hermitian matrices, whose solution involves free probability
theory.  
\head{2.  Large matrices}\endhead
 A frequent  question, occuring both in 
  mathematics and  physics, is to determine the spectrum
of the sum
of  two hermitian matrices. Knowing the respective spectra of two 
$N\times N$ hermitian
matrices $A$
and $B$, the set of
possible spectra for $A+B$ is a  subset of $\Bbb R^N$, which can be
described explicitly and
depends in a complicated way on the  spectra of $A$ and $B$. This problem is not an
easy one, and indeed
its solution was obtained only quite recently (see e.g. \cite{F} for a
discussion).  When $N$ becomes large, 
however, a remarkable phenomenon occurs, and  it 
turns out that, roughly speaking,
 for almost all choices of the matrices $A$ and $B$, with
given spectra, the spectrum of $A+B$ is essentially the same, and can be
computed explicitly,  without knowing the detailed structure of the matrices
$A$ and $B$ (i.e. their eigenvectors).
 In order to give a rigourous mathematical
statement corresponding to the above assertions,
we shall first   introduce  a probability measure on
the set of
matrices having the same spectrum as $A$.  By the spectral
theorem this  set consists of all
matrices 
$UAU^*$, where $U$ describes the group of unitary $N\times N$ matrices.
It
is the orbit  $\Cal O_A$ of $A$ under the adjoint action of $U(N)$, and as
such it carries a
unique invariant probability measure $\rho_A$, the image of the
 Haar measure on $U(N)$ by the map
$U\mapsto UAU^*$.

Given  an hermitian  matrix 
$A$, with spectrum
$\lambda_1,\ldots,\lambda_N$ (counted with multiplicities),
we denote by $\nu_A$  the empirical distribution
on the set of its
eigenvalues, i.e. the 
probability measure
$$\nu_A={1\over N}\sum_{j=1}^N\delta_{\lambda_j}.$$ 
Knowledge of the spectrum (with multiplicities) and of $\nu_A$ are equivalent.

We can now  state the result we had in mind.
\proclaim{Theorem 1} Let for each $N$ positive integer,  $A_N$ and $B_N$ 
be two hermitian matrices, whose norm is 
 bounded uniformly in $N$, and let $\nu_1$ and $\nu_2$ be two probability
measures
with compact supports on $\Bbb R$, such that   $\nu_{A_N}\to\nu_1$ and
$\nu_{B_N}\to\nu_2$ weakly as $N\to\infty$, then there exists a 
probability measure, depending only on $\nu_1$ and $\nu_2$,
denoted 
$\nu_1\boxplus\nu_2$, such that $\nu_{A'_N+B_N'}\to\nu_1\boxplus\nu_2$
weakly, in
probability, as $N\to\infty$, where $A'_N$ and $B'_N$ are random matrices
chosen
independently, with respective distributions $\rho_{A_N}$ and $\rho_{B_N}$. 
\endproclaim
Thus, if we are given two $N\times N$ hermitian matrices (with $N$ large),
 and we only know 
their spectra, then  we can bet, with a good chance to win,
 that the  measure 
$\nu_{A+B}$ is well approximated by the measure
$\nu_A\boxplus\nu_B$. In other words, loosely speaking,
we can compute the spectrum of $A+B$, knowing only the spectra of $A$ and $B$!

Since any probability measure with compact support can be approximated by
 measures like $\nu_A$ for arbitrarily large matrices, Theorem 1
 introduces a binary operation $\boxplus$ on the set of 
compactly supported
measures on the real line which, for reasons to be explained below, we
shall call the free
convolution of measures. This binary operation is clearly associative and
commutative.  

It is instructive to compare the preceding result with 
the case where all considered matrices are supposed to be diagonal.
 Let $A$ be a self-adjoint diagonal matrix, then the set of
 all diagonal matrices
 having the same spectrum as $A$, the analogue of the set $\Cal O_A$, 
  is an orbit of the symmetric group $S_N$, acting by
 permutation of diagonal entries. Let us denote the normalized
 counting measure on this set by $\xi_A$, then one has, denoting by $\mu*\nu$
 the
 usual convolution of two probability measures on $\Bbb R$,  
 \proclaim{Theorem 2} Let for each $N$,
 $A_N$ and $B_N$ be two real diagonal matrices,
  and let $\nu_1$ and $\nu_2$ be two probability
measures
with compact supports, such that   $\nu_{A_N}\to\nu_1$ and
$\nu_{B_N}\to\nu_2$ weakly as $N\to\infty$, then
$\nu_{A'_N+B_N'}\to\nu_1*\nu_2$
weakly, in
probability, as $N\to\infty$, where $A'_N$ and $B'_N$ are random matrices
chosen
independently, with respective distributions $\xi_{A_N}$ and $\xi_{B_N}$. 
\endproclaim
Although we do not know an adequate reference, this last result is probably
well known. In fact, it is not difficult to deduce it from known
concentration of
measure results on the symmetric group, as e.g. in \cite{M}.

Coming back to free convolution, note that
we have not said   how to compute explicitly the measure
$\nu_1\boxplus\nu_2$ in terms of $\nu_1$ and $\nu_2$. This is where
free probability comes in. As we shall see in the
next section,  
Theorem 1 is the consequence of a more general result concerning
large hermitian
matrices, namely the fact that they give rise, asymptotically, to free
random variables.
\head{3. Free random variables}\endhead
Before we describe the connection with large matrices, we introduce the
necessary notions from 
free probability.
We start with some purely algebraic definitions.

Let $\Cal A$ be a complex  algebra with a unit, and $\varphi$ a $\Bbb C$-valued
linear form
on $\Cal A$,
such that 
$\varphi(1)=1$.  Although, in the examples that we shall consider, the algebra
$\Cal A$ will always
be non-commutative, it will be convenient to think  of the elements of this
algebra as random variables, while the map $\varphi$ should
be considered as the
expectation map of classical probability theory (here and in the sequel, we use
the word classical for usual probability theory, as opposed to the
non-commutative theory we develop).

We now introduce the basic notion of free probability theory.
\proclaim{Definition 1}
Let $I$ be a set of indices, and
$\Cal B_i$, for
$i\in I$, be subalgebras of $\Cal A$, containing the unit, then the algebras
$\Cal B_i;i\in I$ are called 
free if one has
$\varphi(a_1\ldots a_n)=0$ each time $\varphi(a_j)=0$ for all $j=1,\ldots
,n$ and
$a_j\in \Cal B_{i_j}$ for some indices $i_1\not= i_2\not =\ldots \not=i_n$.
\endproclaim

Pursuing our analogy with classical probability, if we think of the
algebras $\Cal B_i$
as algebras of random variables, measurable with respect to some sub-sigma
field, then the
above definition is a non-commutative analogue of the definition of
independent
subalgebras. In fact, although it might not be obvious at first sight, 
this definition 
captures the essence of both the notion  
of algebraic independence, and that of independence of sigma-algebras
in classical probability. Let us stress however, that this definition is
not a non-commutative extension of the notion of independence, indeed 
 algebras generated by independent random variables,
 in the sense of classical probability theory, are not free in the
 sense of the above definition.

Before  making some elementary comments, we shall give a
convenient
definition, whose analogy with classical probability should be obvious.
\proclaim{Definition 2} Let $\Cal H_i;i\in I$ be subsets of $\Cal A$, they
are called free if the
subalgebras $\Cal B_i(=$ unital algebra generated by $\Cal H_i$), for $i\in I$,
 are free.
\endproclaim

In order to get acquainted with freeness we shall make a few computations.
So let $\Cal B,\Cal
C\subset \Cal A$ be two free subalgebras, and $b\in \Cal B$, $c\in \Cal C$,
then we can write
$b=b'+\bar b$ where $\bar b=\varphi(b)1$, so that $\varphi(b')=0$.
Similarly we have 
$c=c'+\bar c$. Then using the definition of freeness, we see that
$\varphi(b'c')=0$, thus 
$$\aligned\varphi(bc)&=\varphi((b'+\bar b)(c'+\bar c))\\
&=\varphi(b'c'+b'\bar c+\bar bc'+\bar b\bar c)\\
&=\varphi(b)\varphi(c)\endaligned$$
Here we have used $\varphi(\bar b c')=
\varphi(\varphi(b)1.c')=\varphi(b)\varphi(
c')=0$, and similarly $\varphi(b'\bar c)=0$.

This means that for two free elements $b$ and $c$, their expectations
factorize, exactly as in
the case of independent variables. Let us now look at some more subtle
example. Take
$b_1,b_2\in \Cal B$, and $c_1, c_2\in \Cal C$, then we can expand
$$\aligned\varphi(b_1c_1b_2c_2)&=\varphi((b_1'+\bar b_1)(c'_1+\bar c_1)(
b_2'+\bar b_2)(c'_2+\bar c_2))\\
&=\varphi(b'_1c'_1b'_2c'_2+b'_1c'_1b'_2\bar c_2+b'_1c'_1\bar
b_2c'_2+b'_1c'_1
\bar b_2\bar c_2+\\  &\qquad
b'_1\bar c_1b'_2c'_2+b'_1\bar c_1b'_2\bar c_2
+b'_1\bar c_1\bar b_2c'_2+b'_1\bar c_1\bar b_2\bar c_2+\\
&\qquad\bar b_1 c'_1b'_2c'_2+\bar b_1 c'_1b'_2\bar c_2+\bar b_1 c'_1\bar
b_2c'_2+
\bar b_1 c'_1\bar b_2\bar c_2+\\
&\qquad\bar b_1 \bar c_1b'_2c'_2+\bar b_1 \bar c_1b'_2\bar c_2
+\bar b_1 \bar c_1\bar b_2c'_2+\bar b_1 \bar c_1\bar b_2\bar
c_2)\endaligned
$$
Using the fact that  $\varphi(b'_1c'_1b'_2c'_2)=0$ by the freeness
property,  we are left with 
terms in which some expectation factorizes. For example we shall
treat the term 
$\varphi(b'_1\bar c_1 b'_2c'_2)=\varphi(c_1)\varphi(b'_1b'_2c'_2)$. By the
factorization property already obtained one has
$\varphi(b'_1b'_2c'_2)=\varphi(b'_1b'_2)\varphi(c'_2)=0$. The other terms can be
treated by similar considerations, and  after some straightforward
manipulations we
arrive at the formula
$$\varphi(b_1c_1b_2c_2)=\varphi(b_1b_2)\varphi(c_1)\varphi(c_2)+
\varphi(b_1)\varphi(b_2)\varphi(c_1c_2)-
\varphi(b_1)\varphi(b_2)\varphi(c_1)\varphi(c_2)$$
We observe that the computation of the expectation of the product
$b_1c_1b_2c_2$ can be
reduced to the computation of expectations in the subalgebras $\Cal B$ and
$\Cal C$. This also shows that the result  is different from the one we
would have obtained with independent (commuting) random variables.

It is not difficult to see that the above computation
 can be generalized. More precisely,
taking care of how the terms are successively reduced, the reader should
check, and it is a good exercise do so, that the following is true
\proclaim{Proposition 1} Let $\Cal B_i;i\in I$ be free subalgebras in $\Cal
A$, 
and $a_1,\ldots, a_n\in \Cal A$ such that for all $j=1,\ldots, n$, one has
$a_j\in\Cal B_{i_j}$
for some
$i_j\in I$. Let $\Pi$ be the  partition of $\{1,\ldots, n\}$ determined by
$j\sim k$ if 
$i_j=i_k$.  For each partition $\pi$ of  $\{1,\ldots, n\}$, let  
$\varphi_{\pi}=\prod_{\{j_1,\ldots ,j_r\}
\in \pi\atop j_1<\ldots <j_r}\varphi(a_{j_1}\ldots a_{j_r})$, then  
there exists  universal coefficients $c(\pi,\Pi)$, such that 
$$\varphi(a_1\ldots a_n)=\sum_{\pi\leq \Pi}c(\pi,\Pi)\varphi_{\pi}$$
where the sum is over partitions $\pi$ which are finer than $\Pi$. 
\endproclaim
In particular, this shows
that 
$\varphi(a_1\ldots a_n)$ can be computed explicitly in terms of the
restriction of $\varphi$ to
the algebras $\Cal B_i$. This is a reminiscence of the fact that the joint
distribution of a family
of independent random variables is completely determined if we know the
distribution of each
of the random variables.

As the example that we treated  above suggests,
 the algorithm we have described for 
computing coefficients $c(\pi,\Pi)$ leads quickly to intractable calculations.
Finding an explicit formula for these coefficients is  a non trivial
combinatorial
problem, which fortunately has been solved  by R. Speicher. We shall come back 
to describe his solution later. Let us just mention for the moment that it
involves a certain class of partitions, called ``non-crossing''.

It is time now to state a result showing that the notion we have
introduced is meaningful,
i.e. that there exist non-trivial examples. In this respect
 the situation is
optimum, and we have
\proclaim{Theorem 3}
Let  $\Cal B_i;i\in I$ be complex unital algebras,
  equipped with normalized (i.e. $\varphi_i(1)=1$)
linear forms 
$\varphi_i$, then there  exists an algebra $\Cal A$, with a normalized
linear form $\varphi$,
and unital injective morphisms $\iota_i:\Cal B_i\to\Cal A$,
satisfying 
$\varphi_i=\varphi\circ \iota_i$, and the algebras
$\iota_i(\Cal B_i);i\in I$,  are free in $\Cal A$.
 \endproclaim

The construction of the algebra $\Cal A$ consists in enforcing the freeness
condition in a rather straightforward way. One
defines  it 
as  the ``free amalgamated product'' of the $\Cal B_i$, more precisely,
denote by $\Cal B'_i$, for all $i\in I$,
 the subspace  of elements in $\Cal B_i$ with 
zero expectation, and then take for $\Cal A$ the direct sum 
of $\Bbb C.1$ ($1$ is the unit in $\Cal A$), and 
all spaces $\Cal B_{i_1}'\otimes \ldots \otimes \Cal B_{i_n}'$, where $i_1\ldots
i_n$ runs over all finite sequences in $I$ satisfying $i_1\ne i_2\ne\ldots\ne
i_n$.
In order to turn $\Cal A$ into an algebra,
we need to specify the product of two elements $(a_1\otimes\ldots \otimes
a_n)(b_1\otimes\ldots\otimes b_m)$. If $a_n$ and $b_1$ belong to distinct 
algebras, then
the product is simply $a_1\otimes\ldots \otimes
a_n\otimes b_1\otimes\ldots\otimes b_m$, which belongs to $\Cal A$.
 If $a_n $ and $b_1$ belong to the same
algebra $\Cal B_i$, then $\varphi_i(a_nb_1)$ is not zero in general, so 
we write $a_nb_1=(a_nb_1)'+\overline{a_nb_1}$, and put
$$\aligned(a_1\otimes\ldots \otimes
a_n)(b_1\otimes\ldots\otimes b_m)=&a_1\otimes\ldots \otimes a_{n-1}\otimes
(a_nb_1)'\otimes b_2\otimes\ldots\otimes b_m
\\&+\varphi_i(a_nb_1)(a_1\otimes\ldots \otimes
a_{n-1})( b_2\otimes\ldots\otimes b_m)\endaligned$$ 
Since the ``length'' of $(a_1\otimes\ldots \otimes
a_{n-1})( b_2\otimes\ldots\otimes b_m)$ is $m+n-2$,
the computation can be repeated, and we are done after a finite number of steps.
 The injection $\iota_i$ is defined by 
$\iota_i(b)=\varphi_i(b).1+b'$, and 
the linear form $\varphi$  by $\varphi(a.1)=a$ for $a\in \Bbb C$, and 
$\varphi=0$
 on each space $\Cal B_{i_1}'\otimes \ldots \otimes \Cal B_{i_n}'$. It is not
difficult to check
that the pair $(\Cal A,\varphi)$ fulfills the required conditions.

\head 4. Free convolution\endhead
We are now going  
to be more specific about the kind of algebras that we are considering.
Since we want to do some probability theory, it will be convenient  to be
able to  say that a
function of a random variable is again a random variable, and that a random
variable 
$X$ has a
distribution, namely a probability measure $\mu_X$, such that
$\varphi(f(X))=\int
f(x)\mu_X(dx)$ for any bounded function $f$. The cleanest way to do that is
to assume that
our algebra $\Cal A$ is a {\sl von Neumann algebra}, which means that it is
an algebra of
bounded
operators in some complex Hilbert space $H$, closed under
taking adjoints of operators,   and under
taking limits for 
the weak operator topology
(i.e. simple weak convergence on $H$).
Since we want $\varphi$ to be an analogue of the expectation map with respect to
a probability measure,
we need some positivity assumption. Positivity here will be taken in the sense
of operator theory, namely an element $X$ of $\Cal A$ will be positive
if it is a self-adjoint positive operator on $H$.
The positivity requirement on $\varphi$ is that it takes nonnegative values on
positive operators in $\Cal A$. We shall also assume that $\varphi$
is continuous for the weak operator topology. 
You do not need to be  familiar with von Neumann
algebras in order
to be able to read the following, in fact the only property of a von
Neumann algebra that we
will use is the stability with respect to  functional calculus, namely if
$X\in
\Cal A$ is a  self-adjoint operator on $H$,  then the operator $f(X)$,
which can be defined by
functional calculus for any bounded Borel function $f$ on $\Bbb R$, still belongs to
the algebra $\Cal A$.
Also, if $X$ is self-adjoint  then the map $f\mapsto \varphi(f(X))$,
defined for Borel
bounded functions, is given by
$f\mapsto\int_{\Bbb R}f(x)\mu_X(dx)$, for a unique probability measure 
$\mu_X$, with compact support
 on $\Bbb R$. The probability measure $\mu_X$ is called the
distribution of $X$. Thus, self-adjoint operators in $\Cal A$ behave somewhat
like (real valued) random variables.

 If $\Cal A$ is a von Neumann algebra, and $\varphi$ is a normalized,
weakly continuous, positive linear form on $\Cal A$, we shall call the couple
$(\Cal A,\varphi)$ a non-commutative probability space.

It turns out that the amalgamated free product of algebras exists
 in the category of 
non-commutative probability spaces, namely if the $(\Cal B_i,\varphi_i)$ in
Theorem 3 are non-commutative probability spaces, then we can choose for
$(\Cal A,\varphi)$ a non-commutative probability space, and the $\iota_i$ are
weakly continuous. This  can be seen  by applying 
 the GNS construction to the amalgamated free product 
 constructed after  Theorem 3. We refer to 
\cite{VDN} for a detailed proof. 

We are now in position to give the free-probabilistic
 definition of the free convolution of measures,
which was introduced in Theorem 1.
Let $\mu$ and $\nu$ be probability measures with compact support on $\Bbb R$,
 then there exists a
non-commutative  probability space $(\Cal A,\varphi)$ and self-adjoint elements
$X$, $Y$  in $\Cal A$ with respective   distributions $\mu$ and $\nu$, such
that  $X$ an $Y$ are free. In order to construct $(\Cal A,\varphi)$ and $X,Y$, we can
do the following:  take the algebra $L^{\infty}(\Bbb R,\mu)$,
which can be considered as a von Neumann algebra of operators on
$L^2(\mu)$,  elements of $L^{\infty}(\Bbb R,\mu)$ acting by multiplication
on $L^2(\mu)$. The expectation map $\varphi_{\mu}$ (integration with respect to
$\mu$) turns the pair $(L^{\infty}(\Bbb R,\mu),\varphi_{\mu})$ into
a non-commutative
probability space. The map $x\mapsto x$ on $\Bbb R$ defines a self-adjoint
 element of 
$L^{\infty}(\Bbb R,\mu)$ (since $\mu$ has compact support), which we call $X$. 
By the general construction described above, there exists a non-commutative
probability space, containing 
 $(L^{\infty}(\Bbb R,\mu),\varphi_{\mu})$ and 
 $(L^{\infty}(\Bbb R,\nu),\varphi_{\nu})$ as free subalgebras, and thus
  we get a non-commutative probability
 space containing two free elements $X$ and $Y$ with respective distributions
 $\mu$ and $\nu$. Since $X$ and $Y$ are bounded selfadjoint operators,
 the distribution of $X+Y$ is a probability measure with compact support on $\Bbb R$,
  and it
 is characterized by its moments. Expanding
 the expression $\varphi((X+Y)^n)$ and using  Proposition 1, we see that 
 the moments $\varphi((X+Y)^n)$ can be computed as  polynomial functions of the
 moments of $\mu$ and $\nu$, so that the distribution of $X+Y$   
 really depends only on the measures $\mu$ and $\nu$ and not on the particular
 realization of the operators $X$ and $Y$, thus we
can define $\mu\boxplus \nu$ as the distribution of
 $X+Y$, where $X$ and $Y$ is any pair of free random variables, with respective
 distributions
 $\mu$ and $\nu$.
  Clearly this definition makes the operation $\boxplus$ the free analogue
 of the convolution of measures in classical probability.

 We shall now recover  Theorem 1 from the following result, which shows that
 large hermitian
 matrices give an asymptotic model for free random variables.
 \proclaim{Theorem 4}  Let $A_N$,
  $B_N$,  
$\nu_1$ and $\nu_2$ be as in the statement of Theorem 1, and let
$X$, $Y$ be two self-adjoint elements of some
non-commutative probability space $(\Cal
A,\varphi)$, with respective distributions $\nu_1$ and $\nu_2$, then for every
non-commutative polynomial in two indeterminates, $P$, one has
 ${1\over N}tr(P(A_N',B_N'))\to\varphi(P(X,Y))$
in
probability, as $N\to\infty$, where $A'_N$ and $B'_N$ are random matrices
chosen
independently, 
 with respective distributions $\rho_{A_N}$ and $\rho_{B_N}$. 
\endproclaim
Applying  Theorem 4 to the non-commutative polynomials
$P_n(X,Y)=(X+Y)^n$, we see that the moments of
$A_N'+B_N'$ converge in probability towards that of 
$X+Y$. Since we are dealing 
with compactly supported measures, we 
conclude to the weak convergence of the associated distributions,
and we get Theorem 1.

We shall not give a proof of Theorem 4 here, since this is a far from easy
result. This is a consequence of a more general result from the paper
\cite{Vo2} (see also \cite{VDN}),
 and another, more direct proof can also been found
in \cite{X}. 

It is now time we give a formula for computing the free convolution of two
measures. For convolution of probability measures, one way to perform this
computation is to use the Fourier transform, which converts convolution into
multiplication. Observe that for probability measures with compact support,
the logarithm of the Fourier transform can be expanded into a power series
$$\log (\int_{\Bbb R} e^{itx}\mu(dx))=\sum_{n=1}^{\infty}\sigma_n(\mu)(it)^n.$$
This follows from a  substitution of the power series 
$$\int_{\Bbb R} e^{itx}\mu(dx)=
1+\sum_{n=1}^{\infty}{(it)^n\over n!}\int_{\Bbb R}
x^n\mu(dx)$$ into the expansion of $\log(1+z)$.
 The coefficients $\sigma_n(\mu)$ are
polynomial functions in the moments of  $\mu$, called the cumulants of the
measure $\mu$, and the multiplicativity property
of the Fourier transform of convolution shows that these cumulants linearize
the convolution, namely 
$\sigma_n(\mu*\nu)=\sigma_n(\mu)+\sigma_n(\nu)$. It turns out that the free
convolution admits a similar description in terms of so-called ``non-crossing
cumulants'', which are defined as follows. Let 
$$G_{\mu}(\zeta)=\int_{\Bbb R}{1\over \zeta-t}d\mu(t)$$
be  the Cauchy transform of
$\mu$. 
This defines an analytic function on $\Bbb C\setminus \Bbb R$, such that
$G_{\mu}(\bar \zeta)=\overline{ G_\mu( \zeta)}$, and
$G_{\mu}(\Bbb C^+)\subset \Bbb C^-$ where $\Bbb C^+$ and $\Bbb C^-$
denote respectively the sets of  complex numbers with positive and negative
imaginary part. We can expand this function in a Laurent series involving the
moments of $\mu$
$$G_{\mu}(\zeta)=\sum_{n=0}^{\infty}\zeta^{-n-1}\int_{\Bbb R}
x^n\mu(dx)$$
This series can be inverted
 formally with respect to  $\zeta$, the formal inverse having the 
form 
$$K_{\mu}(z)={1\over z}+\sum_{n=1}^{\infty}C_nz^{n-1}$$ where
the coefficients $C_n$ are polynomial functions of the moments $m_k=\int_{\Bbb
R}x^k\mu(dx)$ 
of $\mu$. The most convenient way to proceed in order to compute the $C_n$
is to write the equation defining
$K_{\mu}$ as
$$G_{\mu}(\zeta)K_{\mu}(G_{\mu}(\zeta))=1+\sum_{n=1}^{\infty}C_nG_{\mu}(\zeta)^n
=\zeta G_{\mu}(\zeta)$$ Equating the coefficients of $\zeta^{-n}$
in the second and third member of the above equality,
 we can evaluate
the $C_n$ in terms of the moments recursively.
The first few values are
$$\aligned C_1&=m_1\\
C_2&=m_2-m_1^2\\
C_3&=m_3-3m_1m_2+2m_1^3\endaligned$$
These formulas can be inverted, and we get the moments in terms of the $C_n$ as
$$\aligned m_1&=C_1
\\
m_2=&C_2+C_1^2\\ m_3=&C_3+3C_1C_2+C_1^3\endaligned$$

The coefficients $C_n$, which we shall call the free cumulants of the measure $\mu$ 
play the role of cumulants for the free convolution, namely one has
\proclaim{Theorem 5} For all compactly supported measures
$\mu$ and $\nu$, on $\Bbb R$, and all $n\geq 1$,
one has
$$C_n(\mu\boxplus \nu)=C_n(\mu)+C_n(\nu).$$ 
\endproclaim Since one can recover the moments of a measure 
 from its free cumulants,
this determines
completely the measure $\mu\boxplus \nu$.

This result was first proved by Voiculescu in \cite{Vo1}, using free creation
and annihilation operators (the formula
has also been discovered independently around the same time, 
in the more restrictive
context of random walks on free products of groups, see \cite{W}).
 Later, Speicher used his combinatorial approach
to freeness to give another proof of Theorem 5. We shall describe Speicher's
proof in section 5. 

Let us give an example of the  computation of $\boxplus$.
 Let $$\mu=\nu={1\over
2}(\delta_{0}+\delta_{1})$$ one has
$$G_{\mu}(\zeta)={\zeta-1/2\over \zeta(\zeta-1)}\text{ and } 
K_{\mu}(z)={z+1+\sqrt{1+z^2}\over
2z}$$ This gives  $$K_{\mu\boxplus\mu}(z)=1+{1\over z}\sqrt{1+z^2}$$
 and 
$$G_{\mu\boxplus\mu}(\zeta)={1\over \sqrt{\zeta(\zeta-1)}}$$ 
The Cauchy transform  can be inverted to give
$$\mu\boxplus\mu(dx)={1\over \pi\sqrt{x(2-x)}}dx\text { on } [0,+2].$$
This is the famous arcsine distribution (here on the interval
$[0,+2]$).
If we recall Theorem 3, 
we have the following striking interpretation of this computation:
take a large integer $N$, and two subspaces of $\Bbb C^N$, of dimension
$[N/2]$, then for most choices of these subspaces,
the sum of the corresponding orthogonal projections  has
an eigenvalue distribution which is well approximated by the arcsine
distribution.
\head 4. Theory of addition of free random variables\endhead
Once one has defined the free convolution of measures, one may try to develop
the theory of addition of free random variables in parallel with the 
theory of addition of independent random variables. It turns out that most of
the well known classical results, such as the law of large numbers, the central
limit theorem, or the  L\'evy-Khintchine formula, 
have free analogues, and this develops into a
beautiful new branch of mathematics, which sheds new lights on some well known
results like the Wigner theorem on spectra of random gaussian matrices.
We shall give a brief survey of these results below, but before that we need to
extend the free convolution to probability measures with unbounded
support. 
\subhead 4.1 Free convolution of measures with unbounded support\endsubhead
 Let $\mu$ be an arbitrary probability measure
 on $\Bbb R$ and let  
$$G_{\mu}(\zeta)=\int_{\Bbb R}{1\over \zeta-t}d\mu(t)$$
be its Cauchy transform. Since $\mu$ may very well have no moment at all, we no
longer have the expansion of $G_{\mu}$ into a power series in $\zeta^{-1}$,
however the following is nevertheless true,
let $$\Theta_{\alpha,\beta}=
\{z=x+iy\,\vert \,y<0;\alpha y< x<-\alpha y;\vert z\vert\leq
\beta\}$$
For every $\alpha>0$, there exists a real number
$\beta>0$ such that  the function
$G_{\mu}$ has a right inverse defined on the domain 
$\Theta_{\alpha,\beta}$, taking values in some domain of the form
$$\Gamma_{\gamma,\lambda}=\{z=x+iy\,\vert\, y>0;-\gamma y< x<\gamma
y;\vert z\vert\geq
\lambda\}$$
with $\gamma,\lambda>0 $.  Call
$K_{\mu}$ this right inverse, and let
$R_{\mu}(z)=K_{\mu}(z)-{1\over z}$. The function $R_{\mu}$ is called the
$R$-transform of the measure $\mu$. Given another probability measure
$\nu$ on $\Bbb R$, we shall  define a new probability measure
$\mu\boxplus\nu$ by the requirement that 
$$R_{\mu\boxplus \nu}=R_{\mu}+R_{\nu}$$ on some domain of
the
form
$\Theta_{\alpha,\beta}$, where these three functions are defined.
It turns out that this definition is meaningful, and it is clear that it
coincides with the previous definition in the case where $\mu$ and $\nu$ have
compact support. In fact, there is also an interpretation of $\mu\boxplus\nu$ as
the distribution of the sum of two free (unbounded in general) self-adjoint
elements affiliated to some non-commutative probability space, see \cite{BV} for details.
\subhead 4.2 Law of large numbers\endsubhead
 The most general formulation of the
law of large numbers for sums of independent equidistributed
 random variables asserts that the distribution of ${1\over n}(X_1+\ldots
 +X_n-M_n)$ converges weakly towards the Dirac measure at zero, for some constants
 $M_n$, if and only if the common distribution $\mu$ of the $X_j$ satisfies
 $t\mu(\Bbb R\setminus{[-t,t]})\to 0$ as $t\to\infty$. It turns out that exactly
 the same result holds if one replaces independent random variables by free
 ones. This is proved in \cite{BP}.
 \subhead 4.3 The Central Limit Theorem\endsubhead
  Let $\mu$ be a probability measure having
  zero mean and finite variance $\sigma$, then one has 
  $R_{\mu}(z)\sim \sigma z$ as $z\to 0$. Let now 
   $X_1,\ldots ,X_n, \ldots $ be a sequence of identically distributed
  free random
 variables, with common distribution $\mu$, then it is easy to see
 that the distribution
 of 
 ${X_1+\ldots X_n\over \sqrt n}$ has $R$-transform  
 $\sqrt n R_{\mu}({z\over \sqrt{n}})$, which converges, as $n\to\infty$, towards
 $\sigma z$. Using continuity properties of the $R$-tranform, it is not
 difficult to see that this implies the following free central limit theorem
 \proclaim{Theorem 6} The distribution of
 ${X_1+\ldots +X_n\over \sqrt n}$ converges weakly, as $n\to\infty$, towards the
 distribution with $R$-transform $\sigma z$.\endproclaim
 The distribution appearing in the free central limit theorem can be computed,
 by inverting the $R$-transform, it is the famous Wigner semi-circular
 distribution given by the density
 ${1\over 2\pi \sigma}\sqrt{4\sigma-x^2}$ on the interval
 $[-2\sqrt \sigma, 2\sqrt \sigma]$. This central limit theorem, together with 
 Theorem 3 on asymptotic freeness of random matrices, provides a conceptual 
 framework for the well known result of Wigner, on the asymptotic behaviour of
 the spectra of large random matrices with gaussian entries.
 \subhead 4.4 Infinitely divisible distributions\endsubhead
 There is a notion of infinitely divisible measures for the free additive
convolution, which is the obvious one, namely a probability measure $\mu$
on the real line is said to be freely infinitely divisible if for every
positive integer $n$, there exists a measure $\mu_n$ such that
$\mu_n^{\boxplus n}=\mu$. The following characterization of freely infinitely divisible
measures on $\Bbb R$ has been obtained by Bercovici and Voiculescu \cite{BV}.
\proclaim{Theorem 7}
 A probability measure $\mu$, on $\Bbb R$ is infinitely divisible if and
only if the function $R_{\mu}$ has an analytic continuation to the
whole
of $\Bbb C^+$, with values in $\Bbb C^-\cup\Bbb R$, and one has $$\lim_{y
\rightarrow 0,y\in\Bbb R}yR_{\mu}(
iy)=0$$
Furthermore any analytic function $R$ having the above properties is the
$R$-transform of some probability measure.\endproclaim
Functions such as those appearing in Theorem 7, have a
  Nevanlinna
representation, which can be seen
 here as the free analogue of the L\'evy-Khintchine
formula. More precisely, 
let $R$ be the $R$-transform of  some infinitely divisible
probability measure, then there exists a real number  
$\alpha$, and a finite positive  measure $\nu$,  on $\Bbb R$, such that 
$$R(z)=\alpha+\int_{-\infty}^{+\infty}{z+t\over 1-tz}d\nu(t)$$
 The extreme points in this integral representation have an interpretation
 similar to the one of the classical L\'evy-Khintchine formula, namely
 the function $R(z)=\alpha$ is the $R$-transform of a Dirac mass at the point $\alpha$,
 the function $R(z)=\sigma z$, corresponding to a point mass at zero for $\nu$,
 is the $R$ transform of the semi-circular distribution with variance $\sigma$.
 Finally the probability measure with $R$-transform $R(z)=\lambda{z+t\over 1-tz}$ is
 the free analogue of the Poisson distribution, namely, it is the weak limit, as
 $n\to\infty$ of the measures $((1-{\lambda\over n})\delta_0+{\lambda\over
 n}\delta_t)^{\boxplus { n}}$ (the ``free binomial distributions''). 
 \subhead 4.5 Stable distributions\endsubhead
 One can define stable distributions exactly as for classical convolution,
 namely, a probability measure on $\Bbb R$ is called stable, if and only if
 the set of probability measures obtained from $\mu$ by applying affine
 transformations of $\Bbb R$, 
 is stable under free convolution. The set of all free stable
 probability measures on the real line has been determined by Bercovici and 
 Voiculescu \cite{BV}. It turns out that there is a natural one-to-one
 correspondance between stable distributions and free stable distributions,
 and  the domains of 
 attraction are the same in the classical and the free cases.
In fact,   any free stable distribution is the image by an affine transformation
of a distribution whose $R$-transform belongs to the following list
\roster 
\item
$R(z)=e^{i\pi\theta}z^{\alpha-1}$ where $1<\alpha\leq 2$, and
$(\alpha-2)\leq \theta\leq 0$.
\item $R(z)=a+b\log z$ where $a\in\Bbb C^+\cup\Bbb R$
and $b\geq - \Im a/\pi$.
\item $R(z)=e^{i\pi\theta}z^{\alpha-1}$ where $0<\alpha<1$, and
$1\leq \theta\leq 1+\alpha$.\endroster
  In particular, the Cauchy distribution
  is a free stable distribution of stability index one.
  
  For these results, see \cite{BVB}

\head 5. Speicher's combinatorial approach to freenes\endhead
We shall now describe
 Speicher's combinatorial approach to the computation of
the coefficients $c(\pi,\Pi)$ of Proposition 1, and some applications.
A thorough discussion is given in \cite{Sp2}.

 Let $S$ be a totally ordered set.
A partition of the set $S$ is said to have a crossing   if there exists 
a quadruple $(i,j,k,l)\in S^4$, with $i<j<k<l$, such that $i$
and
$k$  belong to some class of the partition and  $j$ and $l$ 
 belong to another class. If a partition has no crossing, it is called
non-crossing. The set of all non-crossing partitions of
$S$ is denoted by
$NC(S)$, it is a lattice with respect to the dual refinement order
(such that $\pi\leq \sigma$ if $\pi$ is a finer partition than $\sigma$).
\par
When
$S=\{1,\ldots, n\}$, with its natural order, we will use the notation
$NC(n)$. Here is an example with $n=8$,
 $\pi=\{\{1,4,5\},\{2\},\{3\},\{6,8\},\{7\}\}$.
\input epsf.tex
 $$\gather \epsfbox{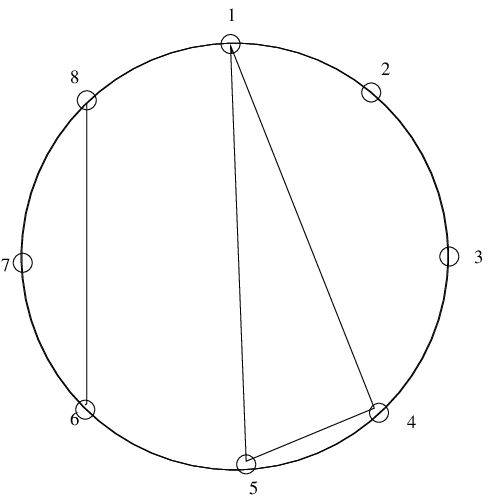}\\ fig.1\endgather$$
 We draw a segment joining each point with the next point in the same 
 class of the partition. The non-crossing condition means that segments should
 not intesect inside the circle.
 
 Let
$(\Cal A,\varphi)$
 be a non-commutative probability space, then  
 we shall define a family  $R^{(n)}$ of
 $n$-multilinear 
 forms on $\Cal A$, for $n\geq 1$,  by the following formula
 $$\varphi(a_1\ldots a_n)=\sum_{\pi\in NC(n)}R[\pi](a_1,\ldots ,a_n)$$
 Here, for $\pi\in NC(n)$, one has 
 $$R[\pi](a_1,\ldots ,a_n)=\prod_{V\in
\pi}R^{(\vert V\vert)}(a_V)$$
where  $a_V=(a_{j_1},\ldots, a_{j_k})$  if $V=\{j_1,\ldots ,j_k\}$ is a
class of the partition $\pi$, with
$j_1<j_2<\ldots <j_k$.  In particular $R[1_n]=R^{(n)}$ if $1_n$ is
the
partition with only one class. 
 For example, one has, for all $a\in \Cal A$, 
 $$\varphi(a)=R(a)$$
  (we forget the superscript $(n)$ when it is clear which $n$
 is considered), for $n=2$
  $$\varphi(a_1a_2)=R(a_1,a_2)+R(a_1)R(a_2)$$ thus
 $$R(a_1,a_2)=\varphi(a_1a_2)-\varphi(a_1)\varphi(a_2)$$
  is the covariance of $a_1$
 and $a_2$. Finally,
 $$\aligned\varphi(a_1a_2a_3)=&R(a_1,a_2,a_3)+
 R(a_1)R(a_2,a_3)+R(a_2)R(a_1,a_3)\\
 &+R(a_3)R(a_1,a_2)+R(a_1)R(a_2)R(a_3)\endaligned$$
 which yields
 $$\aligned R(a_1,a_2,a_3)=&\varphi(a_1a_2a_3)-
 \varphi(a_1)\varphi(a_2a_3)-\varphi(a_2)\varphi(a_1a_3)\\
 &-\varphi(a_3)\varphi(a_1a_2)+2\varphi(a_1)\varphi(a_2)\varphi(a_3)\endaligned
 $$
 In general, for each $n$, one has 
 $$\varphi(a_1\ldots a_n)=R^{(n)}(a_1,\ldots,a_n)+\text{terms involving
 $R^{(k)}$ for $k<n$}$$ so that the $R^{(n)}$ are well
  defined and can be computed by induction on
 $n$.
 
 The non-crossing cumulants can be expressed explicitly in terms of the moments
 by the following formula
 $$R^{(n)}(a_1,\ldots, a_n)=\sum_{\pi\in NC(n)}Moeb(\pi)\varphi
 [\pi](a_1,\ldots ,a_n)$$
 Here $\varphi[\pi](a_1,\ldots,a_n)=
 \prod_{V\in \pi}\varphi(a_{j_1}\ldots a_{j_k})$ where
  $V=\{j_1,\ldots ,j_k\}$  are the classes of $\pi$, and $Moeb$ is the M\"obius
  function on $NC(n)$ defined by $$
  Moeb(\pi)=\prod_{V\in \pi}(-1)^{\vert V\vert}
  c_{\vert V\vert-1}$$
  where $c_n={(2n)!\over n!(n +1)!}$ is the $n^{th}$ Catalan number. 
 
  The connection
between non-crossing cumulants and freeness is the following result of
section 4 of \cite{Sp1}.
\proclaim{Proposition 2}Let $(\Cal B_i)_{i\in
I}$ be free subalgebras of
$\Cal A$, and $a_1,\ldots,
 a_n\in \Cal A$ be such that  $a_j$ belongs to some $\Cal B_{i_j}$
for each
$j\in\{1,2,\ldots ,n\}$. Then $R(a_1,\ldots, a_n)=0$ if there exists some
$j$
and
$k$ with
$i_j\not=i_k$.\endproclaim

Proposition 2 is the key to the computation of the coefficients $c(\pi,\Pi)$ of
Proposition 1. Indeed, let $a_1,\ldots,a_n$ be an arbitrary sequence in $\Cal
A$, 
such that each $a_j$ belongs to one of the algebras $\Cal B_i$, then we have
$$\varphi(a_1\ldots a_n)=\sum_{\pi\in NC(n)}R[\pi](a_1,\ldots ,a_n)$$
and
in this sum the terms corresponding to partitions
$\pi$ having a class containing two elements $j,k$ such that $a_j$ and $a_k$ 
belong to distinct algebras give a zero contribution. Thus we have to sum over
partitions in which all $j$ belonging to a certain class are such that $a_j$
belongs to the same algebra, and the value of $R[\pi](a_1,\ldots,a_n)$
can be expanded in terms of the restriction of $\varphi$ to each of the
subalgebras. 

The multilinear forms $R^{(n)}$ allow us to recover the free cumulants, indeed
one has 
$$\int_{\Bbb R}x^n\,\mu(dx)=\varphi(X^n)=\sum _{\pi\in NC(n)}
 R[\pi](X,\ldots ,X) $$
 \proclaim{Proposition 3}
 Let $X\in\Cal
A$ be selfadjoint and  have distribution $\mu(dx)$, then
 the free cumulants of the measure $\mu$ are 
 given by the formula $C_{n}(\mu)=R^{(n)}(X,\ldots,X)$, 
 for $n=1,2\ldots$.\endproclaim
  Using Propositions 2 and 3, we can now give a proof of Theorem 5.
  Let $X$ and $Y$ be free random variables with
  respective  distributions $\mu$ and $\nu$,
  then the cumulants of $\mu\boxplus\nu$ are given, according to Proposition 3,
   by 
  $$C_n(\mu\boxplus\nu)=R^{(n)}(X+Y,\ldots, X+Y)$$
  Since $R^{(n)}$ is an $n$-linear form, we can expand 
  $R^{(n)}(X+Y,\ldots, X+Y)$, into a sum of 
   terms $R^{(n)}(Z_1,Z_2, \ldots,Z_n)$, where each $Z_i$ is either $X$ or $Y$.
   Applying Proposition 2, we see that all these terms are zero, except
   $R^{(n)}(X,\ldots, X)$ and $R^{(n)}(Y,\ldots, Y)$, thus we have
   $$R^{(n)}(X+Y,\ldots, X+Y)=R^{(n)}(X,\ldots, X)
   +R^{(n)}(Y,\ldots, Y)$$
   and Theorem 5 follows from Proposition 3 again.
   The proofs of Propositions 2 and 3 can be found in \cite{Sp1} or \cite{Sp2}.
   
   \head {6. Some  further topics}\endhead 
   \subhead 6.1 Multiplicative free convolution\endsubhead
   Given two unitary elements $U,V$,
    which are free in some non-commutative probability
   space $(\Cal A,\varphi)$, we can form their product, which is again a unitary
   element. The distributions of $U$ and $V$ are this time probability measures
   say $\mu$ and $\nu$, 
   on the set $ T$ of complex numbers of modulus one, and the distribution of $UV$,
   which depends only on $\mu$ and $\nu$, can be computed by the so-called
   $\Sigma$-transform. Let us introduce the 
$\psi$-function of a probability measure $\mu$, on $ T$,  by
$$\psi_{\mu}(z)=\int_{ T}{z\xi\over 1-z\xi}d{\mu}(\xi)$$
This is a convergent power series in  $D=\{z\in\Bbb C\,\vert
\,\vert z\vert <1\}$, the open unit disk of 
$\Bbb C$, such that $\psi_{\mu}(0)=0$. Let $\Cal M_*$ be the set of probability
measures on $ T$ such that $\int_{ T}\xi d\mu(\xi)\not=0$. If
$\mu\in\Cal M_*$, the function ${\psi_{\mu}\over 1+\psi_{\mu}}$ has a right
inverse, called 
$\tilde{\chi}_{\mu}$, defined in a neighbourhood of $0$,  such that
$\tilde \chi_{\mu}(0)=0$, and we let
$\Sigma_{\mu}(z)= {1\over z}\tilde\chi_{\mu}(z)$ be the
$\Sigma$-transform of $\mu$. Then, for any measures
$\mu,\nu\in\Cal M_*$,  one has $\mu\boxtimes \nu\in\Cal M_*$ and 
$$\Sigma_{\mu\boxtimes\nu}(z)= \Sigma_{\mu}(z)\Sigma_{\nu}(z)$$
in some neighbourhood of zero where these three functions are defined. 

This formula was first found by Voiculescu in \cite{V2}. His proof was
quite complicated, and a  simpler one has been given 
 by U. Haagerup. A proof using
non-crossing partitions, due to  Nica and Speicher is in \cite{NS1}.

There is an analogue, for free multiplicative convolution on $ T$, of the
L\'evy-Khinchine formula. This states that a probability measure on $ T$
is infinitely divisible, for the free multiplicative convolution,  if and
only
if its $\Sigma$ transform can be written as 
$\Sigma_{\mu}(z)=\exp(u(z))$ where $u$ is an analytic function on $D$,
taking values with nonnegative real parts. Such a function has a  
representation of the form 
$$u(z)=i\alpha +\int_{ T}{1+\zeta z\over 1-\zeta z}d\nu (\zeta)$$
for some finite positive measure $\nu$ on $ T$, and real number
$\alpha$.

Finally one can also define multiplicative free convolution for measures on
$\Bbb R_+$.
We shall refer to \cite{BV} for these topics.

We shall here  make a remark on some features which distinguish free
probability from classical probability. Let $\mu$ be a probability 
measure on $\Bbb R$, and let $\Cal G$ be (a suitable branch of)
the logarithm of its Fourier transform, then the set of positive real numbers
$t$ such that $t\Cal G$ is again the logarithm of the Fourier transform of some
probability measure is a closed additive subsemigroup of    $\Bbb R_+$,
containing the positive integers. It is
equal to $\Bbb R_+$, exactly when $\mu$ is infinitely divisible, but it can
 also be reduced to the set of positive integers (e.g if $\mu$ is a Bernoulli
 distribution). In free probability, the analogue of the function $\Cal G$ is
 the $R$-transform. Again, the set of $t$ such that $tR_{\mu}$ is the
 $R$-transform of some probability measure is a closed additive subsemigroup of 
 $\Bbb R_+$, but this time  this subsemigroup always contains the interval
 $[1,+\infty[$. In fact the measure with $R$-transform $tR_{\mu}$ has a nice
 description in terms a free compression. Namely, let $X$ be self-adjoint, with 
 distribution $\mu$ in $(\Cal A,\varphi)$, and let $\pi$ be a self-adjoint
 projection in $\Cal A$, free with $\Cal A$, and such that $\varphi(\pi)={1\over
 t}$, with $t\in [1,\infty[$. The distribution of $\mu$ is the Bernoulli
 distribution $(1-{1\over t})\delta_0+{1\over t}\delta_1$. The set 
 $\pi\Cal A\pi$ of elements of the form $\pi Z\pi$, for $Z\in \Cal A$, is
 an algebra, with $\pi$ as a unit, and  $(\pi\Cal A\pi, t\varphi)$ is a 
 non-commutative probability space. The distribution of 
 the element
 $t\pi X\pi$ in $(\pi\Cal A\pi, t\varphi)$, has a distribution whose
 $R$-transform is given by $tR_{\mu}$ (try to show this using what you know
 about $R$-transforms and non-crossing partitions!). This shows that 
 $tR_{\mu}$ is the $R$-transform of some probability measure, for all $t\geq 1$.
 As the example of Bernoulli distribution shows, there are probability measures
 $\mu$ for which $tR_{\mu}$ is the $R$-transform of some probability measure
 only for $t\in\{0\}\cup[1,\infty[$.
 
 In order to close this section, let us note that the problem of finding the
 distribution of the commutator of two free self-adjoint random variables has
 been solved recently by Nica and Speicher \cite{NS2}.
 
\subhead 6.2 More about random matrices\endsubhead
As we have seen in section 2, free probability provides us with a good
understanding of the way that spectra of large matrices behave under addition.
So far we have not said anything about eigenvectors. It turns out that free
probability again has something to tell us about this. We shall again consider
two large hermitian matrices $A$ and $B$, whose spectra are known. Thinking of
$A+B$ as a perturbation of the matrix $A$, we would like to know how the
eigenvectors of $A+B$ are related to those of $A$. It is clear that if $B$ is
small compared to $A$, then the eigenvectors of $A+B$ should be close to those 
of $A$. Let us denote by $\gamma_1,\ldots, \gamma_N$ the normalized 
eigenvectors of 
$A$, associated with the eigenvalues $\lambda_1,\ldots, \lambda_N$, whereas we
denote by $\xi_1,\ldots, \xi_N$ the eigenvectors of $A+B$, associated with the
eigenvalue $\nu_1,\ldots,\nu_N$. The passage from the
old basis to the new basis is given by the transition matrix $\langle
\xi_k,\gamma_l\rangle$. In fact, since the eigenvectors are only defined up to
some complex number of modulus one, we shall only consider the numbers
$\vert \langle
\xi_k,\gamma_l\rangle\vert^2$, which  form a
bistochastic matrix. These numbers are quite hard to tackle,
 and they may not have
a definite asymptotic behaviour, as $N\to\infty$, so we shall evaluate them
again some test functions. Let $f$ and $g$ be smooth functions on $\Bbb R$, then
we shall look at the asymptotic behaviour of the expression
$$\sum _{1\leq k,l\leq N}g(\nu_k)f(\lambda_l)\vert \langle
\xi_k,\gamma_l\rangle\vert^2$$
We can rewrite this as 
$$tr(g(A+B)f(A))$$ then using our result on asymptotic behaviour of large
matrices, it is easy to see that, if the empirical distributions on the
eigenvalues of $A$ and $B$ converge, as $N\to\infty$, towards
$\mu$ and $\nu$, and we choose matrices at
random as in Theorem 3, then the expression 
$${1\over N}tr(g(A'_N+B'_N)f(A'_N))$$
will converge, in probability, as $N\to\infty$, towards 
$$\varphi(g(X+Y)f(X))$$
where $X$ and $Y$ are free self-adjoint elements, with respective distributions  
$\nu_1$ and $\nu_2$. In principle, we can compute the value of such an
expression, for example if $f$ and $g$ are polynomials. In fact it is easy, 
by a
positivity argument, to see that this value is given by 
$\int f(x)g(u)\rho(dx,du)$ where $\rho$ is some probability measure on $\Bbb
R^2$. It turns out that 
this probability measure can be desintegrated along the values of $x$ as
$\rho(dx,du)= k(x,du)\mu(dx)$, where $k(x,du)$ is a Markov transition kernel,
which could be thought of as the limit of the bistochastic matrices
$\vert \langle
\xi_k,\gamma_l\rangle\vert^2$, and 
this Markov kernel can be explicitly computed,
namely 
it is characterized by  

$$\int_{\Bbb R}(\zeta-u)^{-1}k(x,du)=(F(\zeta)-x)^{-1}\qquad 
\text{for all }\zeta\in \Bbb
C\setminus \Bbb R$$
for some 
function F  on 
$\Bbb C\setminus \Bbb R$ such that  
 $F(\bar\zeta)=\overline{F(\zeta)}$, \ \  $F(\Bbb C^+)\subset \Bbb C^+$,
$Im(F(\zeta))\geq Im(\zeta)$  for $\zeta\in\Bbb C^+$, and
${F(iy)\over iy}\rightarrow 1$ as $y\rightarrow+\infty$, $y\in\Bbb R$,
and
$$G_{\mu}( F(\zeta))=G_{\mu\boxplus \nu}(\zeta)$$
for all $\zeta\in\Bbb C\setminus \Bbb R$.
The map $F$ is uniquely determined by these properties. 
   
   This result appears in \cite{Bi}, where it is applied to the theory
   of processes with free increments.
\refstyle{A}


\widestnumber \key{VDN}
\Refs\nofrills{References}
\ref\key Bi\by P. Biane\paper Processes with free increments
\jour Math. Z.\year 1998\vol 227\pages 143--174\yr 1998\endref
\ref \key BP\by H. Bercovici and V. Pata\jour Ann.  Proba.
\paper The law of large numbers for free identically distributeds
random variables\vol 24 \yr 1996\pages 453--465
\endref
\ref \key BPB\by H. Bercovici and V. Pata\jour Ann. Math. \toappear
\paper Domains of attraction of free stable distributions; with an appendix by
P. Biane on the density of free stable distributions\vol \yr 
\pages \endref
\ref \key BV\by H. Bercovici and D. Voiculescu\jour Indiana University 
Mathematics Journal
\paper Free convolution of measures with unbounded support\vol 42\yr 1993
\pages 733--773\endref
\ref \key F\by W. Fulton\paper Eigenvalues of sums of Hermitian matrices
(after A. Klyachko)\inbook S\'eminaire Bourbaki, expos\'e 845,
Juin 1998\endref
\ref \key M\by B. Maurey\paper Construction de suites sym\'etriques\jour 
C. R. Acad.
Sc. S\'er. A \vol 288\pages 679-681\yr 1979\endref
\ref\key NS1\by A. Nica and 
R. Speicher\paper On the multiplication of free n-tuples of non-commutative
random variables\jour Amer. J. Math.\yr 1996\vol 118\pages 799--837\endref
\ref\key NS2\by A. Nica and 
R. Speicher\paper Commutators of free  
random variables\vol 92\pages 553--592
\jour Duke  Math. J.\yr 1998\endref
\ref\key Sp1\by R. Speicher
\paper Multiplicative functions on the lattice of
non-crossing partitions and free convolution
\jour Math. Annalen\yr 1994\vol 298\pages 141--159\endref
\ref\key Sp2\by R. Speicher\book Combinatorial theory of the free product with
amalgamation and operator-valued 
free probability theory\bookinfo Mem. A.M.S.
\yr 1998\vol 627\endref
\ref\key Vo1\by D. V. Voiculescu\paper Addition of non-commuting random
variables\jour J. Operator Theory\vol 18\yr 1987\pages 223--235\endref
\ref\key Vo2\by D. V.  Voiculescu\paper Limit laws for random matrices and free
products\yr 1991\jour Invent. Math.\vol 104\pages 201--220\endref
\ref\key VDN\by  D. V. Voiculescu,  K. Dykema
and  A. Nica\book Free random variables
\bookinfo  CRM Monograph Series No. 1\publ Amer. Math. Soc.\publaddr
Providence, RI\yr 1992\endref
\ref\key W\by W. Woess\paper 
Random walks on infinite graphs and groups -- asurvey on selected topics
\jour Bull. London Math. Soc.\vol 26\pages 1--60\yr 1994\endref
\ref \key X\by F. Xu\paper A random matrix model from two-dimensional
Yang-Mills theory.\jour Comm. Math. Phys.\yr 1997\vol 190\pages 287--307\endref
\endRefs
\enddocument